\documentclass[a4paper]{article}
\usepackage{amsmath,amssymb,amsfonts}
\usepackage{color}

\newtheorem{theorem}{Theorem}[section]

\title{The antinomy of the liar and provability}
\author{Jailton C. Ferreira}

\date{ }
\begin{document}
\maketitle \pagenumbering{arabic}

\begin{abstract}
This work evidences that a sentence cannot be denominated by $p$ and written
as $p$ \textit{is not true}. It demonstrates that in a system in which $q$ denominates the sentence $q$ \textit{is not provable} it is not provable that $q$ is true and not provable. 
\end{abstract}

\section{Introduction} \label{sec-1}

\hspace{22pt} The antinomy of the liar is an argument regarding the sentence of the liar. The oldest statement known and the most familiar statement of the sentence of the liar are:

\begin{center}
\textit{A man says that he is lying. What he says is true or false?}
\end{center}

and

\begin{center}
\textit{This sentence is not true}
\end{center}

At least since the fourth century B.C., the antinomy of the liar has been discussed without conclusive results.

\hspace{22pt} There are several proposals to treat the antinomy of the liar; among which are the following ones: (i) the liar's antinomy does not proceed, because the liar's sentence does not make sense, for instance ~\cite{Russell} and ~\cite{Quine}; (ii) the liar's sentence is meaningful, but it is neither true nor false, for instance ~\cite{Kripke}; (iii) the liar's sentence is meaningful, it is true or else false, but it is incorrect, in the antinomy of the liar, to pass from the falsehood of the sentence to its truth, for instance ~\cite{Prior}; (iv) the sentence of the liar is meaningful, it is true and it is also false, for instance~\cite{Da Costa}, ~\cite{Priest-1}, ~\cite{Priest-2} and ~\cite{Priest-3}; (v) the sentence of the liar is an example that natural language is self-contradictory, for instance ~\cite{Tarski-1}, and (vi) the liar's sentence is  not a proposition, and therefore, no contradiction is obtained, for instance ~\cite{Strawson}.  

\hspace{22pt} The section \textbf{2} evidences that the sentence of the liar cannot be denominated by $p$ and  written as

\begin{center}
\textit{p} is not true
\end{center}

\hspace{22pt} The section \textbf{3} demonstrates that in a system in which is expressible the sentence

\begin{center}
\textit{q} is not provable
\end{center}

denominated by $q$, where the symbol $q$ in the sentence above is the sentence denominated by $q$, it is not provable that $q$ is true and not provable.

\section{On the formulation of the antinomy of the liar} \label{sec-2}

\hspace{22pt} In his article devoted to the problem of the definition of truth, Tarski ~\cite{Tarski-2} limited himself to formalized languages, abandoning the attempt to solve the problem for the natural language. The reason is in the difficulty using the expression \textit{true sentence} in the natural language, in a consistent way with the laws of the logic.

\hspace{22pt} In the analysis of the concept of truth in the natural language, Tarski began with the semantic definition of truth:

\begin{equation}\label{dois-1}
\textit{x } \textrm{is a true sentence if and only if } \textit{p}
\end{equation}

where, for us to obtain particular definitions, we substituted the symbol $p$ by the sentence and we substituted $x$ by the name of the sentence. Tarski used the classic Aristotelian notion of truth to obtain the definition of line \eqref{dois-1}, the classic notion ~\cite{Aristotle} is expressed in the words

\begin{center}
\textrm{To say of what is that it is not, or of what is not
that it is, is false, \\ while to say of what is that it is, or of what is not that it is not, is true.}
\end{center}

We can denote by \textit{names in quotation marks} any sentence name that consists of an expression in quotation marks, where the expression is the object denoted by the name in question. For instance, `it is raining' is the name of the sentence

\begin{center}
\textrm{it is raining}
\end{center}

\hspace{22pt} We adopt the following procedure: the name of the sentence will be put in the left column and the sentence in the right column. Using this presentation form to the example, we have

\begin{equation}\label{dois-2}
\textrm{`it is raining'} \qquad \qquad \textrm{it is raining}
\end{equation}

For `it is raining' the line \eqref{dois-1} becomes

\begin{equation}\label{dois-3}
\textrm{`it is raining' is a true sentence if and only if it is raining}
\end{equation}

Tarski observes that similar sentences to \eqref{dois-3} seem to be clear and completely in agreement with the meaning of the word truth expressed in \eqref{dois-1}, but that there are situations in which statements of this type in combination with other premises, not less intuitively clear, lead to contradictions. As an example, Tarski presents the following formulation of the antinomy of the liar ~\cite{Tarski-1}, page 158:

\begin{quote}
Let be the sequence of characters \textcolor{red}{++ ==}.

Let be the symbol $c$ the typographic abbreviation for  `the sentence on the fifth non-blank line counting from the first occurrence of the sequence of characters \textcolor{red}{++ ==} in this paper'. Let us consider the following sentence
\begin{center}
\textit{c } \textrm{is not a true sentence}
\end{center}

Considering the meaning of the symbol $c$, we can establish empirically:

\begin{equation}\label{dois-4}
\textrm{`} \textit{c } \textrm{is not a true sentence' is identical to } \textit{c}
\end{equation}

For a name in quotation marks of sentence $c$ (or for any another of their names) we give an explanation of the type \eqref{dois-1}:

\begin{equation}\label{dois-5}
\textrm{`} \textit{c } \textrm{is not a true sentence' is true if and only if } \textit{c } \textrm{is not a true sentence}
\end{equation}

Combining the premises \eqref{dois-4} and \eqref{dois-5}, we obtain the contradiction:

\begin{equation}\label{dois-6}
\textit{c } \textrm{is a true sentence if and only if } \textit{c } \textrm{is not a true sentence}
\end{equation}

\end{quote}

Tarski affirms that (i) the source of the contradiction above was, while building \eqref{dois-5}, the substitution of the symbol $p$ in \eqref{dois-1} for an expression which in itself contains the term ``true sentence'' and that (ii) no rational ground can be given so that such substitutions should be forbidden in principle.

\hspace{22pt} To build four sentences afterwards, we will use the typographic abbreviation `in section \textbf{2} of this article' for `in section \textbf{2} of the article \textit{The antinomy of the liar and provability} found in ArXiv, on June 2008.

\hspace{22pt} Let be

\begin{center}
\textcolor{red}{the sentence in red in section \textbf{2} of this article is not true}
\end{center}

Apparently the sentence immediately above can be denominated `the sentence in red in section \textbf{2} of this article' and rewritten as

\begin{center}
\textrm{`the sentence in red in section \textbf{2} of this article' is not true}
\end{center}

Now let be

\begin{equation}\label{dois-9}
\textit{s} \qquad \qquad 	\textrm{the sentence in red in section \textbf{2} of this article is not true}
\end{equation}

The sentence denominated $s$ affirms something about an indicated sentence. Examining section \textbf{2} of this article we verified that there is a single sentence in red in the section, in other words, the indication made characterizes itself as the name of a sentence. The sentence that $s$ denominates can be rewritten as

\begin{equation}\label{dois-11}
\textit{s} \qquad \qquad 	\textrm{`the sentence in red in section \textbf{2} of this article' is not true}
\end{equation}

Substituting in \eqref{dois-1} the symbol $x$ for $s$ and the symbol $p$ for the sentence denominated by $s$, we have 

\begin{equation}\label{dois-12}
\textit{s } \textrm{is true if and only if `the sentence in red in section \textbf{2} of this article' is not true}
\end{equation}

\eqref{dois-12} means that $s$ is not equivalent to `the sentence in red in section \textbf{2} of this article'. The name `the sentence in red in section \textbf{2} of this article' does not substitute name $s$ nor can be substituted by it. Preserving this distinction is fundamental.

\section{Self-referential sentences and provability} \label{sec-3}

\hspace{22pt} We examined the denominated sentence `the sentence in red in section \textbf{2} of this article' in section \textbf{2} and we concluded the sentence denominated by $s$ is not `the sentence in red in section \textbf{2} of this article'. Now we will examine what happens when the self-referential sentence is of the form

\begin{equation}\label{tres-1}
\textit{p} \qquad \qquad 	\textit{p } \textrm{is not } \alpha
\end{equation}

where $\alpha$ is different from \textit{true} or \textit{not true} or from something equivalent to either \textit{true} or \textit{not true}. Let us substitute recursively $p$ in the sentence for the sentence denominated by $p$:

\hspace{141pt} $p$	\qquad \qquad $p$ is not $\alpha$

\hspace{141pt} $p$	\qquad \qquad ($p$ is not $\alpha$) is not $\alpha$

\hspace{141pt} $p$	\qquad \qquad (($p$ is not $\alpha$) is not $\alpha$) is not $\alpha$

\begin{equation}\label{tres-2}
\hspace{50pt} ....................................................................
\end{equation}

Let us consider sentence $p$ in \eqref{tres-2} with the representations abbreviated for

\begin{equation}\label{tres-3}
\textit{p} \qquad \qquad 	\textit{B } \textrm{is not } \alpha
\end{equation}

If $p$ is true, then it is true that $B$ is not $\alpha$. Since all forms of $B$ have the same truth value and one of the forms is $p$, we have that $B$ is true. If $p$ is false, then we have that $B$ is false.

\hspace{22pt} Now Let us substitute $\alpha$ by ``provable'', that is,

\begin{equation}\label{tres-4}
\textit{p} \qquad \qquad 	\textit{p } \textrm{is not provable}
\end{equation}

In this case the forms of $B$ are

\hspace{160pt} $p$

\hspace{160pt} ($p$ is not provable)

\hspace{160pt} (($p$ is not provable) is not provable)

\begin{equation}\label{tres-5}
\hspace{50pt} ....................................................................
\end{equation}

From \eqref{tres-5} we have

\begin{equation}\label{tres-6}
\textit{p } \textrm{is true if and only if ((} \textit{p } \textrm{is not provable) is not provable) is true}
\end{equation}

The argument above demonstrates the following theorem:

\begin{theorem} \label{theorem-1}
Let be the sentence
\begin{equation}\label{tres-7}
\textit{p } \textrm{is not provable}
\end{equation}
denominated by p. In a system in which p is expressible it is not provable that p is not provable if p is true.
\end{theorem}

\hspace{22pt} Likely the sentence \eqref{tres-6} we have that

\begin{equation}\label{tres-8}
\textit{p } \textrm{is true if and only if (} \textit{p } \textrm{is not provable) is true}
\end{equation}

is part of the truth on \eqref{tres-4}.

\end{document}